\documentclass[twocolumn]{autart}    

\pdfoutput=1
\usepackage{graphicx}          
\newcommand{\der}{{\rm d}}
\usepackage{amsmath}
\usepackage{amsfonts}
\usepackage{euler}
\usepackage{amssymb}
\usepackage{epsfig}
\usepackage{rotating}
\usepackage{graphics}
\newtheorem{theorem}{Theorem}[section]
\newtheorem{lemma}[theorem]{Lemma}
\newtheorem{proposition}[theorem]{Proposition}

\theoremstyle{definition}
\newtheorem{definition}[theorem]{Definition}

\theoremstyle{remark}
\newtheorem{remark}[theorem]{Remark}
\newcommand{\bma}{\begin{pmatrix}}
\newcommand{\ema}{\end{pmatrix}}
\usepackage{pict2e}
\usepackage{calrsfs}
\DeclareMathAlphabet{\pazocal}{OMS}{zplm}{m}{n}
\newcommand{\Da}{\mathcal{D}}
\newcommand{\Db}{\pazocal{D}}

\makeatletter
\DeclareRobustCommand{\loplus}{\mathbin{\mathpalette\dog@lsemi{+}}}
\DeclareRobustCommand{\lotimes}{\mathbin{\mathpalette\dog@lsemi{\times}}}
\DeclareRobustCommand{\roplus}{\mathbin{\mathpalette\dog@rsemi{+}}}
\DeclareRobustCommand{\rotimes}{\mathbin{\mathpalette\dog@rsemi{\times}}}

\newcommand{\dog@rsemi}[2]{\dog@semi{#1}{#2}{-90,90}}
\newcommand{\dog@lsemi}[2]{\dog@semi{#1}{#2}{270,90}}
\newcommand{\dog@semi}[3]{%
  \begingroup
  \sbox\z@{$\m@th#1#2$}%
  \setlength{\unitlength}{\dimexpr\ht\z@+\dp\z@\relax}%
  \makebox[\wd\z@]{\raisebox{-\dp\z@}{%
    \begin{picture}(1,1)
    \linethickness{\variable@rule{#1}}
    \roundcap
    \put(0.5,0.5){\makebox(0,0){\raisebox{\dp\z@}{$\m@th#1#2$}}}
    \put(0.5,0.5){\arc[#3]{0.5}}
    \end{picture}%
  }}%
  \endgroup
}
\newcommand{\variable@rule}[1]{%
  \fontdimen8  
  \ifx#1\displaystyle\textfont3\else
    \ifx#1\textstyle\textfont3\else
      \ifx#1\scriptstyle\scriptfont3\else
        \scriptscriptfont3\relax
  \fi\fi\fi
}
\makeatother
\newcommand{\hook}{\raisebox{-0.45ex}{\makebox[0.6em][r]
{\scriptsize $-$}}\hspace{-0.15em}\raisebox{0.25ex}{\makebox[0.4em][l]{\tiny
 $|$}}}
\newcommand{\om}{\omega}
\newcommand{\be}{\begin{equation}}
\newcommand{\ee}{\end{equation}}
\newcommand{\Span}{\mathrm{Span}}
\newcommand{\dz}{\wedge}
\newcommand{\bbS}{\mathbb{S}}
\newcommand{\bbT}{\mathbb{T}}
\newcommand{\bbR}{\mathbb{R}}
\newcommand{\bbC}{\mathbb{C}}

\newcommand{\sog}{\mathbf{SO}}

\newcommand{\glg}{\mathbf{GL}}
\newcommand{\slg}{\mathbf{SL}}

\newcommand{\gla}{\frak{gl}}
\newcommand{\sla}{\frak{sl}}

\begin{document}

\begin{frontmatter}

\title{Ants and bracket generating distributions \\in dimension 5 and 6} 

\thanks[footnoteinfo]{Support: The research leading to these results has received funding from the Norwegian Financial Mechanism 2014-2021 with project registration number 2019/34/H/ST1/00636. It was also partially supported by the Polish National Science Centre (NCN) via the grant number 2018/29/B/ST1/02583.}

\author[Trieste]{Andrei Agrachev}\ead{agrachev@sissa.it},    
\author[Warszawa]{Pawe\l~ Nurowski}\ead{nurowski@cft.edu.pl}               

\address[Trieste]{Scuola Internazionale Superiore di Studi Avanzati, via Bonomea, 265 - 34136 Trieste, Italy}  
\address[Warszawa]{Centrum Fizyki Teoretycznej,
Polska Akademia Nauk, Al. Lotnik\'ow 32/46, 02-668 Warszawa, Poland}             


\begin{abstract}     We consider a mechanical system of three ants on the floor, in two situations. In the first situation ants move according to {\bf Rule A}, which forces the velocity of any given ant to always point at a neighboring ant; in the second situation ants move according to {\bf Rule B}, which forces the velocity of every ant to be parallel to the line defined by the two other ants. We observe that {\bf Rule A} equips the 6-dimensional configuration space of the ants with a structure of a homogeneous $(3,6)$ distribution, and that {\bf Rule B} foliates this 6-dimensional configuration space onto 5-dimensional leaves, each of which is equied with a homogeneous $(2,3,5)$ distribution. The symmetry properties and the local invariants of these distributions are determined.

  In the case of {\bf Rule B} we study and determine the singular trajectories (abnormal extremals) of the corresponding distributions. We show that these satisfy an interesting system of two ODEs of Fuchsian type.
\end{abstract}

\end{frontmatter}

\section{Rules of motion}
This article is next in the series of papers \cite{nan,agr1,agr,brya,monty,bln}, \cite{hn,en2,en1} presenting simple nonholonomic mechanical systems which are homogeneous models of various parabolic geometries \cite{Cap}. The recent examples of such systems included in this series are very good tools to view in physical terms the main concepts of the parabolic geometry theory. Rather than model the dynamics of the mechanical systems in question, in our case instead of describing the dynamics of the movements of the ants, these papers provide visual examples to illustrate geometry of nonholonomic constraints, and
make a direct connection between an abstract mathematical theory (parabolic geometry) and physics (nonholonomic mechanics).

Although nonholonomic mechanical systems are only briefly mentioned in the usual classical mechanics university courses there are plenty of them in real life. Even those with linear nonholonomic constraints are in abundance: such systems like a skate blade on the plane, a car, trailers, robotic joints and many other man created devices provide examples. Also modeling of movements of animals, such as a movement of a snake, or a falling cat, results in studying systems with nonholonomic constraints.

This short note provides yet another set of examples of this sort. More specifically, we consider three \emph{trained} ants on the floor, idealized as three points $\vec{r}_1$, $\vec{r}_2$ and $\vec{r}_3$ on the plane, which move according to the rules imposed on them by their trainer.

The rules are the {\bf Rules A} and {\bf B} below, and we will analize two \emph{separate} situations: that the ants move by obeying either {\bf rule A} \emph{only}, or {\bf rule B} \emph{only}. Here are the rules:

\begin{itemize}
\item[]{\bf Rule A:} At every moment of time the velocity vector of a given ant, $\frac{\der \vec{r}_i}{\der t}$, should be \emph{aligned} with the direction $\vec{r}_{i+1}-\vec{r}_i$ of the line defined by the ant at $\vec{r}_i$ and the next ant at $\vec{r}_{i+1}$.
\item[]{\bf Rule B:} At every moment of time the velocity vector of a given ant, $\frac{\der \vec{r}_i}{\der t}$, should be \emph{parallel} to the direction $\vec{r}_{i+1}-\vec{r}_{i+2}$ of a line defined by the other two ants.
\item[] $\bullet$ In both rules, $i,j=1,2,3$, and \emph{the sum of two indices}, $i+j$ is \emph{counted modulo} 3. We will assume this convention about sums of indices also in the following.
\end{itemize}

Regardless whether the situation is governed by {\bf rule A} or {\bf B}, the configuration space $M$ of the considered mechanical system is \emph{six} dimensional. It can be, for example, (locally) parametrized by six real numbers corresponding to the $2\times 3=6$ coordinates $(x_1,y_1,x_2,y_2,x_3,y_3)$ of the three points $\vec{r}_i=(x_i,y_i)$ in a chosen Cartesian coordinate system $(x,y)$ on the plane. In this parametrization the movement of the system of ants is described in terms of a curve $m(t)=(x_1(t),y_1(t),x_2(t),y_2(t),x_3(t),y_3(t))$, and its velocity at time $t$ is given by $\dot{m}(t)=(\dot{x}_1(t),\dot{y}_1(t),\dot{x}_2(t),\dot{y}_2(t),$ $\dot{x}_3(t),\dot{y}_3(t))$.

Now, since {\bf rule A} imposes that $\frac{\der\vec{r}_i}{\der t} || (\vec{r}_{i+1}-\vec{r}_i)$ and {\bf rule B} imposes that $\frac{\der\vec{r}_i}{\der t} || (\vec{r}_{i+1}-\vec{r}_{i+2})$, we see that the movement of ants under the {\bf rule A} has velocities constrained according to:
  $$(y_{i+1}-y_i)\dot{x}_i-(x_{i+1}-x_i)\dot{y}_i=0,\,\, i=1,2,3,\quad\quad {\bf Rule A},$$
    and that the movement of ants under the {\bf rule B} have velocities constrained according to:
    $$(y_{i+1}-y_{i+2})\dot{x}_i-(x_{i+1}-x_{i+2})\dot{y}_i=0,\,\, i=1,2,3,\,\, {\bf Rule B}.$$
    In both cases the velocity constraints of the systems, the \emph{nonholonomic} constraints as they are called, are \emph{linear}. Thus the space of admisible velocities at each point $q$ of the configuration space $M$ is a \emph{vector subspace} $D_q$ of the tangent space $\mathrm{T}_qM$. Since in both cases we have \emph{three} independent velocity constraints at each point $q\in M$ the vector spaces $D_q$ are 3-diemensional, and as such, collected point by point, define rank \emph{three} distributions $\Db$ on $M$.

    Let us first make a brief analysis of the geometry of the pair $(M,{\Db})$ in the case of ants moving under {\bf rule A}.
\section{The rule `every ant is chased by precisely one other ant' results in a $(3,6)$ distribution}

    In case of {\bf rule A} the distribution $\Db$ of admissible velocities on $M$ is given by the \emph{annihilator} of the following three 1-forms:
    $$\om_i=(y_{i+1}-y_i)\der x_i-(x_{i+1}-x_i)\der y_i, \quad i=1,2,3,$$
    or which is the same, is spanned by the three vector fileds
    \be Z_i=(x_{i+1}-x_i)\partial_{x_i}+(y_{i+1}-y_i)\partial_{y_i},\quad i=1,2,3\label{zi}\ee
    on $M$,
    $${\Db}=\Span(Z_1,Z_2,Z_3).$$
    Taking the commutators (Lie brackets) of the vector fields $Z_1,Z_2,Z_3$ spanning the distribution $\Db$ we get three new vector fields
    $$\begin{aligned}
    &Z_{i,i+1}=[Z_i,Z_{i+1}]=\\&(x_{i+1}-x_{i+2})\partial_{x_i}+(y_{i+1}-y_{i+2})\partial_{y_i},\quad i=1,2,3.\end{aligned}
    $$
    Now, calculating $Z_1\dz Z_2\dz Z_3\dz Z_{12}\dz Z_{31}\dz Z_{23}$, one gets
  $$\resizebox{0.46\textwidth}{!}{$\begin{aligned}&Z_1\dz Z_2\dz Z_3\dz Z_{12}\dz Z_{31}\dz Z_{23}=\\&\Big(\sum_{i=1}^3(y_ix_{i+1}-x_iy_{i+1})\Big)^3\partial_{x_1}\dz\partial_{y_1}\dz\partial_{x_2}\dz\partial_{y_2}\dz\partial_{x_3}\dz\partial_{y_3},\end{aligned}$}$$
    so it follows that the \emph{six} vector fields $Z_1,Z_2,Z_3,Z_{12},Z_{31},$ $Z_{23}$ are \emph{linearly independent} at each point $m$ of the configuration space $M$, \emph{except} the points on the singular locus, where coordinates of $m$ satisfy
    \be 32~A=\sum_{i=1}^3(y_ix_{i+1}-x_iy_{i+1})=0.\label{area}\ee
    Since the number $A$ defined above is the \emph{area} of the triangle having the three ants as its vertices, we see that this happens for those configurations when three ants stay on a line.

    At this stage we recall a concept of the \emph{growth vector} of a vector distribution. Given any vector distribution $\Db$ we consider the following sequence $\{\Db^I\}$  of derived distributions: The 0-th element of this sequence is $\Db^0=\Db$, and distributions $\Db^I$, with $I=0,1,\dots,$ are defined recursively as $\Db^{I+1}=[\Db^I,\Db^I]+\Db^I$. Here $[\Db^I,\Db^I]$ denotes the set consisting of all linear combinations\footnote{with functional coefficients} of all commutators of vector fields forming the distribution $\Db^I$. The \emph{growth vector} of $\Db$ is a vector $(r_0,r_1,\dots,r_K)$, whose integer components $r_I$ are the respective \emph{ranks} of the derived distributions $\Db^I$, $r_I=\mathrm{rank}(\Db^I)$. Note that the growth vector typically vary from point to point, but there are important examples of distributions, such as those considered in this paper, when the growth vector is \emph{constant} over large open sets. Note also that the last number $r_K$ in the sequence $(r_0,r_1,\dots,r_K)$ can not excceed the dimension $n$ of the manifold $M$ on which the distribution $\Db$ resides, $r_K\leq n$. If the last component $r_K$ of the growth vector of a distribution $\Db$ is equal to the dimension $n$ of the manifold $M$, then the distribution $\Db$ is called \emph{bracket generating}. If a bracket generating distribution $\Db$ is in addition a velocity distribution of a mechanical system, such system is \emph{controllable}.

    Having said this, we summerize our observations about the velocity distribution $\Db$ of the three ants moving under {\bf rule A}: The \emph{growth vector} of this distribution is $(3,6)$ everywhere, except those points in the configuration space which correspond to the three ants staying on a line. A short term for describing this property of $\Db$ is to say that $\Db$ is a $(3,6)$ distribution. In particular, the control system
    $\dot{q}=u_1Z_1(q)+u_2Z_2(q)+u_3Z_3(q)$,
   associated with the ants moving under {\bf rule A}, is controllable.

Rank 3 distributions have \emph{differential invariants} \cite{bry}. We recall, that two distributions ${\Db}_1$ and ${\Db}_2$ on respective manifolds $M_1$ and $M_2$ are (locally) equivalent\footnote{The equivalence of distributions $\Db_1$ and $\Db_2$ is denoted by $\phi_* {\Db}_1={\Db}_2$ to remind that the equivalence is obtained by the diffeomorphism $\phi$. The notation encodes the property that the vectors spaninng the distribution $\Db_1$ at point $q$ in $M_1$, after being transformed by $\phi$ to the point $\phi(q)$ in $M_2$, span there the same vector subspace of the tangent space $\mathrm{T}_ {\phi(q)}M_2$ as the vector space $\Db_2(\phi(q))$ of the distribution $\Db_2$.}, if and only if there exists a (local) diffeomorphism $\phi:M_1\to M_2$ transforming distribution $\Db_1$ to distribution $\Db_2$. In particular the statement about rank 3 distributions having invariants, means that there are locally \emph{non}equivalent rank 3 distributions on 6-dimensional manifolds. Actually, there are \emph{infinitely many} of nonequivalent ones. Among them the $(3,6)$ distributions are generic, and the growth vector $(3,6)$ distinguishes them locally from, for example, distributions with the growth vector $(3,5)$; these latter distributions are rank 3 distributions $\Db$ in dimension 6 such that in the sequence ${\Db}^0={\Db}$, ${\Db}^{I+1}=[{\Db}^I,{\Db}^I]+{\Db}^I$, with $I=0,1,....$, the distribution ${\Db}^1$ is \emph{involutive}\footnote{Some mathematicians prefer the term \emph{integrable} here.} and has rank 5. More importantly, there are locally nonequivalent $(3,6)$ distributions.

    One way of characterizing distributions locally is to determine their Lie algebra of symmetries. Given a manifold $M$ and distribution $\Db$, the Lie algebra of symmetries of $\Db$ consists of vector fields $X$ on $M$ such that $[X,{\Db}]\subset\Db$. Here $[X,\Db]$ denotes the space consisting of commutators (Lie brackets) of the vector field $X$ with all vector fields belonging to $\Db$. Equivalently, a vector field $X$ is a symmetry of $\Db$, if the local flow $\phi_t^X$ of $X$ satisfies $(\phi_t^X)_*\Db\subset\Db$. It is known \cite{brythesis,bry}  that for rank 3 distributions with the growth vector $(3,6)$ the maximal algebra of symmetries is attained for the distribution locally given in Cartesian coordinates $(q^i,p_j)$ in $\bbR^6$ as the annihilator of three 1-forms $\lambda_i=\der p_i+\epsilon_{ijk} q^j\der q^k$, $i=1,2,3$. Here we used the Einstein summation convention stating that repeated indices are summed over; we also used the totally skew-symmetric Levi-Civita symbol $\epsilon_{ijk}$ in $\bbR^3$, which is zero when any two indices are the same, and otherwise, is +1 or -1, depending on the sign of permutation of the three indices ${ijk}$. This maximally symmetric distribution has its Lie algebra of symmetries isomorphic to the 21-dimensional Lie algebra $\mathfrak{spin}(4,3)$ (see \cite{brythesis,bry}).

    Since the velocity distribution $\Db$ of the system of three ants moving according {\bf rule A} has growth vector $(3,6)$ almost everywhere, it is interesting to ask what is its Lie agebra of symmetries. THe answer is given by the following Theorem.
    \begin{theorem}\label{thm21}
    The Lie algebra of \emph{all} symmetries of the velocity distribution $\Db$ of the system of three ants moving according {\bf rule A} is isomorphic to the Lie algebra $\mathfrak{sl}(3,\bbR)$. In coordinates $(x_1,y_1,x_2,y_2,x_3,y_3)$ in $\bbR^6$, as in \eqref{zi}, the 8 independent local symmetries of ${\Db}=\Span(Z_1,Z_2,Z_3)$ are:
    $$\resizebox{0.47\textwidth}{!}{$\begin{aligned}
 X_1=&\partial_{x_1}+\partial_{x_2}+\partial_{x_3},\\
 X_2=&\partial_{y_1}+\partial_{y_2}+\partial_{y_3},\\
 X_3=&y_1\partial_{x_1}+y_2\partial_{x_2}+y_3\partial_{x_3},\\
 X_4=&x_1\partial_{y_1}+x_2\partial_{y_2}+x_3\partial_{y_3},\\
 X_5=&x_1\partial_{x_1}+x_2\partial_{x_2}+x_3\partial_{x_3},\\
      X_6=&y_1\partial_{y_1}+y_2\partial_{y_2}+y_3\partial_{y_3},\\
      X_7=&x_1y_1\partial_{x_1}+x_2y_2\partial_{x_2}+x_3y_3\partial_{x_3}+y_1^2\partial_{y_1}+y_2^2\partial_{y_2}+y_3^2\partial_{y_3},\\
      X_8=&x_1^2\partial_{x_1}+x_2^2\partial_{x_2}+x_3^2\partial_{x_3}+x_1y_1\partial_{y_1}+x_2y_2\partial_{y_2}+x_3y_3\partial_{y_3}.
      \end{aligned}$}$$
    \end{theorem}
\begin{remark}
    We proved this theorem \emph{by explicitly solving the symmetry equations} $[X,\Db]\subset\Db$ for the velocity distribution $\Db$ on $\bbR^6$, with coordinates $(x_1,y_1,x_2,y_2,x_3,y_3)$, as in \eqref{zi}. A reader can convince himself that the eight vector fields $X_i$ from the Theorem are really the symmetries of $\Db$ by taking Lie brakets of $X_i$s with all the velicity distribution generators $Z_i$ and observing that all these brackets belong to $\Db=Span(Z_1,Z_2,Z_3)$. Another way of seeing that $\sla(3,\bbR)$ is included in the symmetry algebra is to observe that the mechanical system of ants moving by {\bf rule A} is defined in terms of no other notions than \emph{points}, \emph{lines} and their \emph{incidence relations}\footnote{such as `a point lying on a line' or `lines intersecting themselves or not', etc} in the plane. These are notions of the \emph{projective geometry} on the plane. A convenient \emph{model} of this geometry is the geometry of \emph{points and lines} lying in the plane $z=1$ of the Cartesian space $\bbR^3$ with coordinates $(x,y,z)$. In this model \emph{lines} and \emph{points} in the $z=1$ plane are the respective intersections of those planes or lines in $\bbR^3$, which contain the $\bbR^3$ space origin $(x,y,z)=(0,0,0)$. Since the Lie group $\glg(3,\bbR)$ naturally acts on the planes and lines passing through the origin in $\bbR^3$, it also acts on the intersection of these planes and lines with the plane $z=1$, namely on the \emph{lines} and \emph{points} in the $z=1$ plane. This $\glg(3,\bbR)$ action on the lines and points lying in the $z=1$ plane is \emph{not effective}, as the scaling group element, a multiple $\lambda I$ of the identity in $\glg(3,\bbR)$, stretches planes and lines in $\bbR^3$ merely, and in turn does not move the corresponding lines and points in the $z=1$ plane. The group that acts effectively on lines and points in the $z=1$ plane, is a group ${\bf PGL}(3,\bbR)$, which is an 8-dimensional \emph{quotient} Lie group $\glg(3,\bbR)/(\lambda I)$. And this group, the \emph{projective linear group} ${\bf PGL}(3,\bbR)$, is the \emph{symmetry group} of the \emph{projective geometry on the plane}, i.e. of the \emph{geometry of lines and points in the plane}. Therefore if we look for a group of symmetries of the ants-under-{\bf rule A}-distribution-$\Db$, a distribution that is entirely defined in terms of projectively invariant notions such as \emph{points}, \emph{lines} and \emph{their incidence in the plane}, this group must contain ${\bf PGL}(3,\bbR)$ - a group whose \emph{Lie algebra} is the Lie algebra $\sla(3,\bbR)$. The local symmetries of the distribution are given in terms of the Lie algebra of its group of (local) symmetries. Therefore the ants' distribution $\Db$ must be \emph{at least} $\sla(3,\bbR)$ \emph{symmetric}. Actually, our Theorem \emph{says more}: the Lie algebra of symmetries of $\Db$ is \emph{equal} to $\sla(3,\bbR)$. This does not follow from the simple arguments mentioned in this remark. We proved it by explicitely solving the symmetry equations.
\end{remark}

\begin{remark}
  Thus although the symmetry of this $(3,6)$ distribution is far from being maximal among all $(3,6)$ distributions, the ants distribution $\Db$, considered in this section, can be locally identified with one of the \emph{homogeneous models} of $(3,6)$ distributions, a model that lives on the homogeneous manifold $PGL(3,\bbR)/\bbT^2$, where $\bbT^2$ is the maximal torus in $PGL(3,\bbR)$. That, this is the case is obvious from the explicit formulas for the symmetry vector fields $X_i$ in Theorem \ref{thm21}. Indded, the symmetry vector field $X_1$ is a local version of the \emph{scaling of all coordinates} $x_1$, $x_2$ \emph{and} $x_3$ \emph{by the same number}; likewise  the symmetry vector field $X_2$ is a local version of the \emph{scaling of all coordinates} $y_1$, $y_2$ \emph{and} $y_3$ \emph{by (perhaps) another number}. Looking at the equations \ref{zi} defining generators $Z_i$ of the distribution, one sees that such \emph{scallings} do not change these generators at all. This shows that these two symmetries belong to the \emph{isotropy subalgebra}, and since $[X_1,X_2]=0$, the corresponding local \emph{isotropy group} is a direct sum $H\oplus H$ of two copies of a 1-dimensional Lie group $H$ \footnote{All 1-dimensional Lie groups are locally isomorphic!}. Thus, if we want to have a \emph{global model} of the ants-under-{\bf rule A}-manifold-$M$, i.e. a 6-dimensional manifold $\tilde{M}$ whose local portions are in one to one correspondence with those parts of $M$ on which the distribution $\Db$ is $(3,6)$, we may take $\tilde{M}={\bf PGL}(3,\bbR)/\bbT^2$. This is because $H\oplus H$ globalizes to the 2-dimensional torus $\bbT^2=\bbS^1\oplus \bbS^1$, where the circle $\bbS^1$ has e.g. a group structure of the \emph{complex numbers} $z\in\bbC$, \emph{whose modulus is equal to one}, $|z|=1$.
\end{remark}

\begin{remark}\label{syn}
We close this section, with a remark that the vector space over the real numbers spanned by the symmetry vector fields $X_1,X_2,X_3,X_4,X_5-X_6$ form a Lie algebra isomorphic to the \emph{semidirect product} of the simple Lie algebra $\sla(2,\bbR)$ and the commutative Lie algebra $\bbR^2$. This product is denoted by $\sla(2,\bbR)\roplus\bbR^2$, so that
$$\sla(2,\bbR)\roplus\bbR^2=\Span_\bbR(X_1,X_2,X_3,X_4,X_5-X_6).$$
Here vector fields $X_1$ and $X_2$ on $\bbR^6$ correspond to translations in the plane in respective directions $\partial_x$ and $\partial_y$. The vector fields $X_3$, $X_4$ and $X_5-X_6$ correspond to the linear transformations of the plane with unit determinant. In particular we have the following identifications of the respective Lie algebra elements: $$X_3\sim\bma 0&1\\0&0\ema,\,\, X_4\sim\bma 0&0\\1&0\ema\,\mathrm{and}\,\, X_5-X_6\sim\bma 1&0\\0&-1\ema.$$
\end{remark}

\section{The rule `each ant moves in a parallel to the line defined by the other two' is not so simple}\label{B}
Now, applying {\bf rule B} to the movement of the three ants, we find that their velocity distribution $\Db$ is given by the annihilator of the three 1-forms
$$\om_i=(y_{i+1}-y_{i+2})\der x_i-(x_{i+1}-x_{i+2})\der y_i, \quad i=1,2,3.$$
It can be spanned by the three vector fileds
    \be Z_i=(x_{i+1}-x_{i+2})\partial_{x_i}+(y_{i+1}-y_{i+2})\partial_{y_i},\quad i=1,2,3\label{zia}\ee
    on $M$,
    $${\Db}=\Span(Z_1,Z_2,Z_3).$$
    The commutators of the vector fields $Z_1,Z_2,Z_3$ spanning $\Db$ are
    \be
    \begin{aligned}
    Z_{i,i+1}&=[Z_i,Z_{i+1}]=\\&(x_i-x_{i+2})\partial_{x_i}+(x_{i+2}-x_{i+1})\partial_{x_{i+1}}+\\&(y_i-y_{i+2})\partial_{y_i}+(y_{i+2}-y_{i+1})\partial_{y_{i+1}},\quad i=1,2,3.
    \label{zij}\end{aligned}\ee
   And now the story is \emph{different} than in the case of {\bf rule A}. Calculating $Z_1\dz Z_2\dz Z_3\dz Z_{12}\dz Z_{31}\dz Z_{23}$, one gets
  $$ Z_1\dz Z_2\dz Z_3\dz Z_{12}\dz Z_{31}\dz Z_{23}=0.$$
   So the rank of the derived distribution ${\Db}^1=[{\Db},{\Db}]+\Db$ is \emph{smaller} than 6. The velocity distribution $\Db$ for the {rule B} is \emph{not} bracket generating! Actually one easilly finds that there is \emph{precisely} one linear relation between the vector fields $(Z_1,Z_2,Z_2,Z_{12},Z_{31},Z_{23})$, namely
   \be Z_1+Z_2+Z_3+Z_{12}+Z_{31}+Z_{23}=0.\label{rel}\ee
   This shows that the velocity distribution $\Db$ for the ants moving under {\bf rule B} has the growth vector $(3,5)$. The first derived distribution ${\Db}^1$ has rank 5 and is \emph{involutive}! The 6-dimensional configuration space $M$ of ants being in a motion obeying {\bf rule B} is foliated by 5-dimensional leaves. Once ants are in the configuration belonging to a given 5-dimensional leaf in $M$ they can \emph{not} leave this leaf by moving according {\bf rule B}! This system is \emph{non}-controllable: its configuration space splits into 5-dimensional orbits.

   Now the question arises about the function that enumerates the leaves of the foliation of the distribution ${\Db}^1$. What is the feature of motion of the ants whose preservation forces the ants to stay on a given leaf?

   There is a quick algebraic answer to this question:

   To see it, note that
   $$\der(\om_1+\om_2+\om_3)=0.$$
   this means that that there exists a function $F$ such that
   $$\der F=\om_1+\om_2+\om_3.$$
   One can directly check that
   $$F~=~32~A,$$
   where $A$ is as in \eqref{area}. Since all three vector fields $Z_i$ as in \eqref{zia} annihilate $\om_i$, and thus they annihilate the 1-form $\om_1+\om_2+\om_3=32\der A$, and in turn they annihilate the one form $\der A$, then they \emph{are tangent} to the 5-dimensional submanifolds $A=const$ in $M$.

   This shows that the ants under {\bf rule B} move in a way such that \emph{the triangle having them as its vertices has always the same area}! This proves the folowing proposition.
   \begin{proposition}
     The triangle with vertices formed by three ants moving according {\bf rule B} has in every moment of time the same area.
   \end{proposition}
   Apart of the algebraic proof of this proposition given above, it can be also seen by a `pure thought' observing that any movement of the three ants obeying {\bf rule B} is a superposition of three primitive moves: an ant $\# i$ moves, and ants $\# (i+1)$ and $\# (i+2)$ rest, for each $i=1,2,3$. In each of the three primitive situations, since the vertex $\# i$ of the triangle moves in a line parallel to the corresponding base $\# (i+1)-\# (i+2)$ of the triangle, the area of the triangle formed by the ants $\# 1$, $\# 2$ and $\# 3$ is obviously unchanged. Since the general movement according to {\bf rule B} is a linear combination of the three primitive movements preserving the area, it also preserves the area.

   So we see that the movement of the ants according to {\bf rule B} stratifies the configuration space: once in an initial position the ants defined a triangle $\Delta$ of area $A$, they move on a 5-dimensional submanifold $M_A$ of $M$ whose configuration points correspond to triangles $\Delta'$ having the same area $A$ as $\Delta$. For each fixed $A$, the three vector fields $(Z_1,Z_2,Z_3)$ as in \eqref{zia} are \emph{tangent} to the five manifold $M_A$. They define a distribution ${\Db}=\Span(Z_1,Z_2,Z_3)$ there, whose growth vector is $(3,5)$.

   The 3-distribution $\Db$ on each leaf $M_A$ is actually a \emph{square} of a rank 2-distribution $\Da$. By this we mean that there is a rank 2-distribution $\Da$ such that its first derived distribution $\Da^1=[\Da,\Da]+\Da$ equals $\Db$.
   Indeed, consider
   $$\Da=\Span(Z_1-Z_2,Z_3-Z_1)$$
   with $Z_1,Z_2,Z_3$ as in \eqref{zia}. Since $[Z_1-Z_2,Z_3-Z_1]=-Z_{12}-Z_{31}-Z_{23}$, with $Z_{ij}$ as in \eqref{zij}, then using relation \eqref{rel} we get $$[Z_1-Z_2,Z_3-Z_1]=Z_1+Z_2+Z_3$$ and consequently $$[Z_1-Z_2,Z_1-Z_3]\dz (Z_1-Z_2)\dz (Z_3-Z_1)=3Z_3\dz Z_2\dz Z_1.$$ This shows (i) that for each $A=\mathrm{const}$ the commutator $[\Da,\Da]$ is tangent to $M_A$ and (ii) that the first derived distribution of $\Da$ on $M_A$ is the entire 3-distribution, $[\Da,\Da]+\Da=\Db$.
   Thus we have just established the following proposition.
   \begin{proposition}
     The 6-dimensional configuration space $M$ of three ants moving on the plane according to {\bf rule B} is foliated by 5-dimensional submanifolds $M_A$ consisting of configuration points defining triangles of equal area $A$ on the plane. The ants obeying {\bf rule B} must stay on a given leaf $M_A$ of the foliation during their motion. Their velocity distribution $\Db$ of rank 3, defines a rank 2 distribution $\Da$, which is the `square root' of $\Db$, $$\Db=[\Da,\Da]+\Da.$$ The rank 2 distribution $\Da$ has the growth vector $(2,3,5)$ on each leaf $M_A$.
     \end{proposition}
   We recall that rank 2 distributions with growth vector $(2,3,5)$ on 5-dimensional manifolds have local differential invariants \cite{cartan}. In particular their symmetry algebra can be as large as 14-dimensional Lie algebra $\mathfrak{g}_2^*$ of the split real form\footnote{Similarly to the complex Lie group $\sog(3,\bbC)$ which has two real forms $\sog(3)$ and $\sog(1,2)$, the complex simple exceptional Lie group ${\bf G_2}$ has also \emph{two} real forms. One of these forms is \emph{compact}, like $SO(3)$, and the other is non-compact, like $SO(1,2)$. This second non-compact real form of complex group ${\bf G_2}$ is denoted by ${\bf G_2^*}$, with the Lie algebra denoted by $\mathfrak{g}_2^*$; in the structural theory of simple Lie groups the real Lie group ${\bf G_2^*}$ is called the \emph{split} real form of complex ${\bf G_2}$. It is this real form of the complex Lie group $G_2$, namely the split $G_2^*$, that is relevant in the geometry of $(2,3,5)$ distributions.} of the simple exceptional complex Lie group ${\bf G_2}$. This happens for the rank 2 distribution given on a 5-dimensional quadric $p_iq^i=1$ in $\bbR^6$, with coordinates $(q^i,p_i)$, as the annihilator of three 1-forms $\lambda_i=\der p_i+\epsilon_{ijk} q^j\der q^k$, $i=1,2,3$.

   The rank 2 distribution $\Da$ on each 5-dimensional leaf $M_A$ has a 5-dimensional Lie algebra of symmetries corresponding to the Lie group of \emph{affine transformations} of the plane \emph{preserving area}. In coordinates $(x_1,y_1,x_2,y_2,x_3,y_3)$ in $M=\bbR^6$ the Lie algebra of these transformations is spanned by the five symmetry vector fields $(X_1,X_2,X_3,X_4,X_5-X_6)$ from Remark \ref{syn}. Denoting by $S$ a vector field
   $$S=a_1 X_1+a_2X_2+a_3 X_3+a_4 X_4+a_5(X_5-X_6),$$
   with $a_\mu=\mathrm{const}$, $\mu=1,2,\dots 5$, one can directly check that for $A$ given by \eqref{area} and for $Z_i$ given by \eqref{zia} we have:
   $$\begin{aligned}S(A)&=0,\\
  & [S,Z_1-Z_2]\dz(Z_1-Z_2)\dz(Z_3-Z_1)=0\\&\hspace{3.cm}\mathrm{and}\\&[S,Z_3-Z_1]\dz(Z_1-Z_2)\dz(Z_3-Z_1)=0.
   \end{aligned}$$
   We invoked this algebraic argument, if the reader would not agree with us that the $\sla(2,\bbR)\roplus \bbR^2$ symmetry of the rank 2 distribution $\Da$ is obvious.

   We are now in a position to state the theorem, which will be proven in the subsequent sections of the article:

   \begin{theorem}\label{33}
     The Lie algebra of \emph{all} symmetries of the velocity distribution $\Da$ of the system of three ants moving on the plane according to {\bf rule B} is isomorphic to the Lie algebra $\sla(2,\bbR)\roplus \bbR^2$ of group of motions on the plane preserving volume. The distribution is one of the homogeneous models of $(2,3,5)$ distribution, which can be locally realized on the 5-manifold being the group $\slg(2,\bbR)\rtimes \bbR^2$.

     The fundamental invariant of $\mathcal D$, its harmonic curvature encapsulated in the so called \emph{Cartan quartic}, is of algebraic type $D$, or what is the same, has no real roots.
   \end{theorem}
   \begin{remark}\label{34re}
     The definition of the Cartan quartic for $(2,3,5)$ distributions is beyond the scope of this article (see \cite{cartan}, and \cite{nan2}, p. 94, for details).  Here, we only mention that it gives a neat way of expressing the lowest order local differential invariant of any $(2,3,5)$ distribution. This is given as a certain totally symmetric tensor $C_{ABCD}$, $A,B,C,D=1,2$, with respect to the action of the $\glg(2,\bbR)$ group. The tensor $C_{ABCD}$ has \emph{five} real components $A_1=C_{1111}$, $A_2=C_{1112}$, $A_3=C_{1122}$, $A_4=C_{1222}$, $A_5=C_{2222}$ which are used to encapsulate it in the \emph{quartic}, the \emph{Cartan quartic}, of a $(2,3,5)$ distribution. This quartic reads: $C(z)=A_1+4A_2z+6A_3z^2+4A_4z^3+A_5z^4$, with a complex variable $z\in\bbC$. The \emph{number} and \emph{multiplicities} of \emph{roots} of this quartic, considered as polynomial in variable $z$, provide particularly simple invariants of a $(2,3,5)$ distribution. A number of cases for the number of roots and their multiplicities may happen. Here we only say that the Cartan quartic of a $(2,3,5)$ distribution \emph{is of type} $D$, if the Cartan quartic has \emph{one pair of mutually complex conjugated} roots, \emph{each of them having multiplicity two}.
   \end{remark}
 \begin{remark}
   A reader who is interested in the Cartan's approach to the geometry of the distribution $\mathcal D$ can now jump directly to Section \ref{se5}. Theorem \ref{33} is proven there in the spirit of Cartan's 5-variables paper \cite{cartan}. There is however another approach to the analysis of invariant properties of distributions used by geometric control theorists, which we discuss now. This approach uses the important notion of \emph{singular trajectories} (or \emph{abnormal extremals}). In the next two sections we will determine these trajectories for ants distributions, and will prove Theorem \ref{33} using the ideas related to them.
 \end{remark}
\section{Singular trajectories for ants' movement}

We start this section with generalities about singular trajectories for vector distributions, and then we apply this to determine these trajectories for the distributions associated with the ants movement under the {\bf rules A} and {\bf B}.

Let $\Db\subset TM$ be a smooth vector distribution on a smooth manifold $M$ and $\Db^\perp\subset T^*M$ be its annihilator. We denote by $\Db_0^\perp\subset\Db^\perp$ the bundle $\Db^\perp$ with the removed zero section.

Recall that the cotangent bundle $T^*M$ is equipped with a canonical symplectic form $\sigma$.
\begin{definition} A Lipschitz curve $t\mapsto\lambda(t),\ t\in[0,1],$ in $\Db_0^\perp$ is called a \emph{singular} or an \emph{abnormal extremal} of $\Db$ if $\dot\lambda\in\ker\left(\sigma|_{\Db^\perp}\right)$. The projection of $\lambda(\cdot)$ to $M$ is called singular trajectory or abnormal geodesic.
\end{definition}
Let us explain the geometric meaning of the introduced notions and thus motivate the terminology. We'll do it without going to analytic details which can be found in first chapters of the book \cite{AgBaBo}. Let
$$
\Omega_\Db=\{\gamma:[0,1]\to M \mid \dot\gamma(t)\in\Db_{\gamma(t)}M,\ 0\le t\le 1\},
$$
be the space of integral curves of the distribution $\Db$
and $F_t:\Omega_\Db\to M$ be the evaluation map, $F_t(\gamma)=\gamma(t)$. Then singular trajectories are just critical points of the ``boundary map'' $(F_0,F_1):\Omega_\Db\to M\times M$.

Moreover, let $\gamma$ be a singular trajectory; the curve
$$
t\mapsto\lambda(t)\in T^*_{\gamma(t)}M\setminus\{0\},\quad 0\le t\le 1,
$$
is a singular extremal if and only if
$$
\lambda(t)D_\gamma F_t=\lambda(s)D_\gamma F_s, \quad \forall\,t,s\in[0,1].
$$
Here $D_\gamma F_t:T_\gamma\Omega\to T_{\gamma(t)}M$ is the differential of $F_t$ at $\gamma$,
$\lambda(t):T_{\gamma(t)}M\to\mathbb R$ and $\lambda(t)D_\gamma F_t$ is the composition of the linear map $D_\gamma F_t$
and linear form $\lambda(t)$

We see that $\lambda(t)$ plays the role of ``Lagrange multipliers'' corresponding to critical points. Indeed, $\gamma$ is a critical point of the map $(F_0,F_1)$ if and only if there exists a pair of covectors $(\lambda(0),\lambda(1))\ne 0$ such that $\lambda(0)D_\gamma F_0=\lambda(1)D_\gamma F_1$. Moreover, if $\gamma$ is a critical point of the ``boundary map'', then the restriction of the curve $\gamma$ to any segment $[t,s]\subset [0,1]$ is also a critical point of the boundary map.

To effectively compute singular extremals, we use the Pontryagin Maximum Principle which exploits the Hamiltonian language. Given a smooth function $h:T^*M\to\mathbb R$,
Hamiltonian vector field $\vec h$ on $T^*M$ is defined by the identity $dh=\sigma(\cdot,\vec h)$. The Poisson bracket
 $\{a,b\}$ of two functions on $T^*M$ is defined by the formula
 $
 \{a,b\}=\langle db,\vec a\rangle=\sigma(\vec a,\vec b)
 $
 and provides $C^\infty(T^*M)$ with the structure of a Lie algebra.

 For any subset $S\subset T_\lambda(T^*M)$, we denote by $S^\angle\subset T_\lambda(T^*M)$ the skew-orthogonal complement of $S$,
 $$
 S^\angle=\{\xi\in T_\lambda(T^*M) \mid \sigma(\xi,S)=0\};
 $$
 then $(S^\angle)^\angle=\Span\,S$. If $\lambda$ is a regular point of a function $a:T^*M\to\mathbb R$ and $a(\lambda)=c$,
 then $T_\lambda a^{-1}(c)=\ker d_\lambda a=\vec a(\lambda)^\angle$.

 Assume that vector fields $Z_1,\ldots,Z_k$ on $M$ generate the distribution $\Db$, i.\,e.
 $$
 \Db|_q=\Span\left(Z_1(q),\ldots,Z_k(q)\right),\quad \forall\,q\in M.
 $$
 We define $h_i:T^*M\to\mathbb R,\ i=1,\ldots,k$, by the formula:
 $$
 h_i(\lambda)=\langle\lambda,Z_i(q)\rangle,\quad \forall\,\lambda\in T^*_qM,\ q\in M.
 $$
Then
$$
\Db^\perp=\bigcap\limits_{i=1}^kh_i^{-1}(0),\quad T_\lambda\Db^\perp=\bigcap\limits_{i=1}^k\ker d_\lambda h_i
$$
and
$$
\left(T_\lambda\Db^\perp\right)^\angle=\Span\{\vec h_1(\lambda),\ldots,\vec h_k(\lambda)\}.
$$

We have: $\ker\left(\sigma|_{\Db^\perp}\right)=(T\Db^\perp)\cap(T\Db^\perp)^\angle$. Hence a curve
$t\mapsto\lambda(t)$ in $T^*M$ is a singular extremal if and only if there exist real functions $t\mapsto u_i(t),\ t\in[0,1],$ usually called control function, such that
$$
\dot\lambda(t)=\sum\limits_{i=1}^ku_i(t)\vec h_i(\lambda(t))\quad \mathrm{and}\quad h_i(\lambda(t))\equiv 0,\ i=1,\ldots,k.
$$

More notations: we set $h_{ij}=\{h_i,h_j\},\ H_{IJ}=\{h_{ij}\}_{i,j=1}^k$ and we have:
$$h_{ij}(\lambda)=\langle\lambda,[Z_i,Z_j](q)\rangle,\ \forall\,\lambda\in T^*_qM,\ q\in M.$$

We differentiate identities $h_i(\lambda(t))\equiv 0$ in virtue of the differential equation and obtain
$$\sum\limits_{j=1}^ku_i(t)h_{ij}(\lambda(t))\equiv 0.$$ In other words, $u(t)\in\ker H_{IJ}(\lambda(t))$.

In our models {\bf A} and {\bf B}, $H_{IJ}(\lambda)$ is a nonzero $3\times 3$-matrix, $\forall\,\lambda\in
T^*_qM\setminus \{0\},\ q\in M$. Hence the anti-symmetric matrix $H_{IJ}$ has a one-dimensional kernel spanned by the vector
$(u_1,u_2,u_3)=(h_{23},h_{31},h_{12})$. It follows that any singular extremal, up-to a reparametrization, is a solution of the differential equation
\be
\dot\lambda=h_{23}(\lambda)\vec h_1(\lambda)+h_{31}(\lambda)\vec h_2(\lambda)+h_{12}(\lambda)\vec h_3(\lambda). \label{2}\ee

Moreover, $\Db^\perp$ is an invariant submanifold of the dynamical system \eqref{2} on $T^*M$. Indeed,
$\frac d{dt}h_i(\lambda(t))=h_{il}h_{ji}+h_{ji}h_{li}=0$ in virtue of \eqref{2}, where $i\ne j\ne l$. It follows that, up-to  reparameterizations, singular extremals are exactly trajectories of the dynamical system \eqref{2} starting at nonzero elements of $\Db^\perp$.
Note that the right-hand side of \eqref{2} is anti-symmetric with respect to permutations of the indices $(1,2,3)$.

The projection of equation \eqref{2} to $M$ has a form:
$$
\dot q=h_{23}Z_1(q)+h_{31}Z_2(q)+h_{12}Z_3(q)=\sum\limits_{i=1}^3u_iZ_i(q).
$$
Moreover,
\be\begin{aligned}
\frac d{dt}h_{j,j+1}(\lambda)&=\sum\limits_{i=1}^3u_i\{h_i,h_{j,j+1}\}(\lambda)=\\&
\sum\limits_{i=1}^3u_i\langle\lambda,[Z_i,[Z_j,Z_{j+1}]]\rangle. \label{3}\end{aligned}
\ee

In both models {\bf A} and {\bf B}, \resizebox{0.27\textwidth}{!}{$\Span\left(Z_i, [Z_j,Z_{j+1}],\ i,j=1,2, 3\right)$} is a Lie algebra. This is a 6-dimensional Lie algebra of the group of affine transformations of the plane in the case {\bf A} and a 5-dimensional algebra of the group of area preserving affine transformations of the plane in the case {\bf B}. In the last case we have:
$\sum\limits_{i=1}^3(Z_i+[Z_i,Z_{i+1}])=0$.

Anyway, in both cases vector fields $[Z_i,[Z_j,Z_{j+1}]]$ are linear combinations of the vector fields $Z_1,Z_2,Z_3,$ $[Z_2,Z_3],[Z_3,Z_1],[Z_1,Z_2]$ with constant coefficients. Hence the function $\lambda\mapsto \langle\lambda,[Z_i,[Z_j,Z_{j+1}]]\rangle$ is a linear combination of functions
\be
h_1(\lambda),\ h_2(\lambda),\ h_3(\lambda),\ h_{23}(\lambda),\ h_{31}(\lambda),\ h_{12}(\lambda),\quad \lambda\in T^*M, \label{4}\ee
with the same constant coefficients. If $\lambda\in\Db^{\perp}$, then first three elements of the sequence \eqref{4} vanish and
last three elements are just $u_1,u_2,u_3$. Coming back to \eqref{3} we obtain that $\dot u_i$ is a quadratic function of $u_1,u_2,u_2,\ i=1,2,3$.

A straightforward calculation of Lie brackets gives:
$$
\left\{\begin{aligned}\dot u_1&=-u_1(u_1+u_2)\\ \dot u_2&=-u_2(u_2+u_3)\\ \dot u_3&=-u_3(u_3+u_1),\\ \end{aligned}\right.
\eqno ({\bf A})
$$
$$
\left\{\begin{aligned}\dot u_1&=u_1(u_2-u_3)\\ \dot u_2&=u_2(u_3-u_1)\\ \dot u_3&=u_3(u_1-u_2).\\ \end{aligned}\right.
\eqno ({\bf B})
$$

Now we focus on the case {\bf B} and postpone for the future the apparently more complicated case {\bf A}. Recall that in this case $u_1+u_2+u_3=0$. Moreover, the product $u_1u_2u_3$ is the first integral of system ({\bf B}) and thus the system is integrable. In addition, system ({\bf B}) is anti-symmetric with respect to the permutation of variables and the central reflection $(u_1,u_2,u_3)\mapsto(-u_1,-u_2,-u_3)$. The fact that
$u_1(t)u_2(t)u_3(t)=const$ along trajectories of ({\bf B}) implies that $u_i(t)$ do not change sign. We see that it is enough to study the system in the domain $u_1,u_2\ge 0,\ u_3=-u_1-u_2$.

Let $z_1=x_1+iy_1,z_2=x_2+iy_2,z_3=x_3+iy_3$ be the vertices of the triangle; then $Z_i=\langle z_{i+1}-z_{i+2},\frac\partial{\partial z_i}\rangle,\ i=1,2,3.$ Singular trajectories are solutions of the system:
$$
\left\{\begin{aligned}\dot z_1&=u_1(z_2-z_3)\\ \dot z_2&=u_2(z_3-z_1)\\ \dot z_3&=u_3(z_1-z_2),\\ \end{aligned}\right.
$$
where $u_1,u_2,u_3$ satisfy ({\bf B}).

Let $z_0=\frac 13(z_1+z_2+z_3)$ be the barycenter of the triangle. A direct differentiation gives that $\ddot z_0(t)\equiv 0$, i.e. the barycenter moves along a straight line with constant velocity.

Moreover, given a triangle and a straight line through its barycenter, there is exactly one singular trajectory whose barycenter moves along this line.

We set $\zeta_i=z_i-z_0,\ i=1,2,3,$ so that $\zeta_1+\zeta_2+\zeta_3=0$. Now we eliminate variables $u_3=-u_1-u_2$ and $\zeta_3=-\zeta_1-\zeta_2$ from the differential equations and arrive to the system:
\be
\resizebox{0.43\textwidth}{!}{$\left\{\begin{aligned}\dot u_1&=u_1(u_1+2u_2)\\ \dot u_2&=-u_2(2u_1+u_2),\\ \end{aligned}\right.
\qquad
\left\{\begin{aligned}\dot \zeta_1&=(u_1+u_2)\zeta_1+u_1\zeta_2\\ \dot \zeta_2&=-u_2\zeta_1-(u_1+u_2)\zeta_2.\\ \end{aligned}\right. \label{5}$}\ee
We also have: $\dot z_0=u_1\zeta_2-u_2\zeta_1=const.$ Another polynomial first integral of system \eqref{5} is already mentioned function
$u_1u_2(u_1+u_2)$.

  System \eqref{5} can be easily solved; the solution is expressed in elliptic functions. We denote by $e$ the velocity of the barycenter of the triangle, $e=u_1\zeta_2-u_2\zeta_1$, it does not depend on $t$. We have
$$
\left\{\begin{aligned}\dot \zeta_1&=(u_1+2u_2)\zeta_1+e\\ \dot \zeta_2&=-(2u_1+u_2)\zeta_2+e.\\ \end{aligned}\right.
$$
In other words,
$$
\left\{\begin{aligned}\dot \zeta_1&=\frac{\dot u_1}{u_2}\zeta_1+e\\ \dot \zeta_2&=-\frac{\dot u_2}{u_2}\zeta_2+e.\\ \end{aligned}\right.
$$
It follows that
$$
\begin{aligned}\zeta_1(t)&=\frac{u_1(t)}{u_1(0)}\zeta_1(0)+\int\limits_0^t\frac{u_1(t)}{u_1(\tau)}\,d\tau\, e,\\
\zeta_2(t)&=\frac{u_2(0)}{u_2(t)}\zeta_2(0)+\int\limits_0^t\frac{u_2(\tau)}{u_2(t)}\,d\tau\, e.\\ \end{aligned}
$$

Let $u_1^2u_2+u_1u_2^2=c$; then $2u_1u_2=\sqrt{u_1^2+4u_1c}-u_1^2$ and we obtain
$\dot u_1=\sqrt{u_1(u_1^3+4c)}$. Similarly, $\dot u_2=\sqrt{u_2(u_2^3+4c)}$. Hence
$$
\int\limits_{r_1}^{u_1}\frac{dv}{\sqrt{v(v^3+4c)}}=t=\int\limits^{r_2}_{u_2}\frac{dv}{\sqrt{v(v^3+4c)}},
$$
where $r_1=u_1(0),\ r_2=u_2(0)$. So $t$ is presented as elliptic integrals of $u_1$ and $u_2$.
Relations between the constants:
$$
r_1r_2(r_1+r_2)=c,\quad r_1\zeta_2(0)-r_2\zeta_1(0)=e.
$$

\medskip
Another alternative way to study the same system is to try to eliminate variables $u_1,u_2$ and focus completely on the vertices of the triangles.
Let $\nu_1=u_1+u_2,\ \nu_2=u_1u_2$ be elementary symmetric functions of $u_1,u_2$ and $\delta=u_1-u_2$ the discriminant.
We have:
$$\resizebox{0.48\textwidth}{!}{$\dot\nu_1=\delta\nu_1,\quad \dot\nu_2=-\delta\nu_2,\quad, \dot\delta=\nu_1^2+2\nu_2,\quad \delta^2=\nu_1^2-4\nu_2.$}
$$

The case $\nu_1=0$ corresponds to the constant trajectory (recall that we are working in the domain $u_1,u_2\ge 0$).
Let $\nu_1\ne 0$. The case $\nu_2=0$ corresponds to a singular trajectory with a fixed vertex of the triangle.
This is a singular trajectory whose barycenter moves along the line connecting the barycenter with the vertex.

The equations are very simple in this case. We leave the calculations to a reader as an exercise and formulate only the final result:
Let $z_1(t)=const$, then, for all $t$, the passing through $z_1$ median of the triangle is a segment of one and the same straight
line, the opposite to $z_1$ side of the triangle remains parallel to itself and the length of this side is regulated by the constant area condition.

Now turn to the main case $\nu_2\ne 0$. We see that $\dot\delta>0$; the permutation of $u_1$ and $u_2$ changes the sign
of $\delta$ and we can restrict ourselves to the study of the parts of trajectories in the domain $\delta>0$. Let us make a
time substitution in our system and introduce new time $s$ according to the rule $\frac{ds}{dt}=\delta$. Then
$\frac{d\nu_1}{ds}=\nu_1,\ \frac{d\nu_2}{ds}=-\nu_2$. Hence $\nu_1(s)=c_1e^s,\ \nu_2(s)=c_2e^{-s}$.

Moreover, $u_1=\frac 12(\nu_1(s)+\delta),\ u_2=\frac 12(\nu_1-\delta)$ and we have:
$$
\left\{\begin{aligned}\frac{d\zeta_1}{ds}&=\frac{\nu_1}\delta \zeta_1+\left(\frac{\nu_1}{\delta}+1\right)\frac{\zeta_2}2\\
\frac{d\zeta_2}{ds}&=\left(\frac{\nu_1}{\delta}-1\right)\frac{\zeta_1}2-\frac{\nu_1}\delta \zeta_2.\\ \end{aligned}\right.
$$
Now we make one more change of the time variable:
$$
\tau=\frac{\nu_1(s)}{\delta(s)}=(1-ce^{-3s})^{-\frac 12},
$$
where $c=\frac{4c_2}{c_1^2}$. Then $\frac{d\tau}{ds}=\frac 32\tau(1-\tau^2)$ and we obtain:
$$
\left\{\begin{aligned}\frac{d\zeta_1}{d\tau}&=\frac{2\tau}{3(1-\tau^2)} \zeta_1+\frac{\zeta_2}{3\tau(1-\tau)}\\
\frac{d\zeta_2}{d\tau}&=\frac{\zeta_1}{3\tau(1+\tau)}-\frac{2\zeta_2}{3(1-\tau^2)}\\ \end{aligned}\right.
$$

This is a Fuchsian system with 3 poles $-1,0,1$. We can rewrite it in the canonical matrix form:
$$\resizebox{0.46\textwidth}{!}{$
3 \frac{d\Psi}{d\tau}=\left[\frac 1{\tau-1}\begin{pmatrix}-1& -1\\ 0& 1 \end{pmatrix}+\frac 1\tau\begin{pmatrix} 0& 1\\ 1& 0 \end{pmatrix}+
\frac 1{1+\tau}\begin{pmatrix} 1& 0\\ -1& -1 \end{pmatrix}\right]\Psi,$}
$$
where $\Psi=(\zeta_1,\zeta_2)^*$. Note that the three $2\times2$ matrices appearing in this equation form a basis of the Lie algebra $\sla(2)$.
In principle, this kind of systems (with 3 poles) is resolved in the Gauss hypergeometric function.
The following explicit expression is obtained by Renat Gontsov (Moscow):
$\Psi(\tau)=\frac 1{\sqrt2}\,X(\tfrac{2\tau}{\tau+1})$, where
$$\resizebox{0.49\textwidth}{!}{$
X(t)=(t-1)^{-1/3}\,t^{1/3}\left(\begin{array}{cc}\frac85(t-1)t & 2t-1 \\ & \\
                                                           0 & t-1
                              \end{array}\right)\left(\begin{array}{cc} F(\frac53,\frac73,\frac83;t) & -\frac54\,t^{-5/3} \\ & \\
                                                                        F(\frac23,\frac43,\frac53;t) & t^{-2/3}
                                                      \end{array}\right)$}
$$
and $F(\alpha,\beta,\gamma;t)$ is the Gauss hypergeometric function. Recall that
$$
F(\alpha,\beta,\gamma;t)=\sum_{j=0}^{\infty}\frac{(\alpha)_j\,(\beta)_j}{(\gamma)_j\,j!}\,t^j,
$$
where $(\alpha)_j=\alpha(\alpha+1)\ldots(\alpha+j-1)$. We see that $\gamma-\alpha=1$ in the explicit expression; this corresponds to the so called degenerate case of the hypergeometric equation (see \cite{BaEr} ) and should lead to a further simplification of the expression. It would be nice to get a clear geometric interpretation of the computed singular trajectories.

Actually, all this story is a benefit performance of the number 3! It is continued in the next section.

\section{Sub-Riemannian structure and Cartan quartic}\label{sec4}

Here we show how the structure of the Cartan quartic is determined by the natural symmetries of the distribution $\Db$.

\begin{lemma} Given a triangle $\Delta$ in the plane, there exists a unique ellipse $C_\Delta$ such that $C_\Delta$ contains
the vertices of the triangle and the tangent line to $C_\Delta$ at any vertex is parallel to the opposite side of the triangle.
\end{lemma}
{\bf Proof.} The desired properties are preserved by the affine transformations of the plane. To prove the existence we
transform the triangle in the regular one and then take the circle that contains the vertices of the regular triangle.

To prove the uniqueness we take an ellipse $C$ that satisfies the desired properties and transform it into a circle by an affine transformation. We see that the image of the triangle under this transformation must be regular.
\qquad $\square$

Any ellipse defines an Euclidean metric on the plane: we say that a segment has length 1 if the endpoint of a parallel segment are the center of the ellipse and a point on the ellipse itself.

Now we can define a natural Riemannian structure on the space of triangles that is invariant with respect to the group of affine transformations. Given a curve $t\mapsto\Delta_t=conv\{z_1(t),z_2(t),z_3(t)\}$ in the space of triangles, we define the square of the length of its velocity $(\dot z_1(t),\dot z_2(t),\dot z_3(t))$ as the sum of squares of the lengths of the vectors $\dot z_i(t),\ i=1,2,3,$
measured in the Euclidean structure associated to the ellipse $C_{\Delta_t}$.

The restriction of this Riemannian metric to the distribution $\Db=\Span\{Z_1,Z_2,Z_3\}$ is a sub-Riemannian metric such that the square of the length of the vector $\sum\limits_{i=1}^3u_iZ_i$ is equal to $\sum\limits_{i=1}^3u_i^2$.
This is valid for both models {\bf A} and {\bf B}.

In what concerns the case {\bf B}, velocities of singular trajectories form a rank 2 distribution
${\mathcal D}=\{\sum\limits_{i=1}^3u_iZ_i \mid \sum\limits_{i=1}^3u_i=0\}\subset\Db$. Recall that the vector fields
$\sum\limits_{i=1}^3u_iZ_i,\ \sum\limits_{i=1}^3u_i=0$ generate a 5-dimensional Lie algebra. The famous Cartan invariant of this kind of distribution (Cartan quartic) is a degree 4 homogeneous form on ${\mathcal D}$.

\begin{proposition}\label{pro} The Cartan quartic of ${\mathcal D}$ is equal to the form:
$$
\sum\limits_{i=1}^3u_iZ_i\mapsto c\left(u_1^2+u_2^2+u_3^2\right)^2,\quad \sum\limits_{i=1}^3u_i=0,
$$
where $c$ is a constant.
\end{proposition}

{\bf Proof.} Cyclic permutations of the fields $Z_1,Z_2,Z_3$ induce automorphisms of the Lie algebra $Lie\{Z_1,Z_2,Z_3\}$ and hence symmetries of the distribution ${\mathcal D}$. It follows that the Cartan quartic written as a form of $u_1,u_2,u_3$
must be invariant with respect to the cyclic permutations of the variables $u_1,u_2,u_3$.

Consider a projective line $$\bar{\mathcal D}=\{u_1:u_2:u_3 \mid \sum\limits_{i=1}^3u_i=0\}.$$
The group of cyclic permutations of $u_1,u_2,u_3$ acts freely on $\bar{\mathcal D}$ and preserves the sets of real roots of a prescribed multiplicity of the Cartan quartic. Any orbit of this cyclic group has 3 points, hence the number of roots of a prescribed multiplicity is a multiple of 3. At the same time, if the quartic is not identical zero, then it may have 0,\ 1,\ 2 or 4 roots of a prescribed multiplicity and never 3 roots. It follows that the quartic does not have real roots.

The quartic is real, hence it is a square (or minus square) of a positive definite quadratic form on $\bar{\mathcal D}$.
This quadratic form is invariant with respect to the action of the order 3 cyclic group,hence it is proportional to the form
$\sum\limits_{i=1}^3u_iZ_i\mapsto \sum\limits_{i=1}^3u_i^2$. Indeed a nonempty intersection of the level sets of this form and
any non proportional to it quadratic form has 2 or 4 elements and cannot be preserved by a free action of the order 3 cyclic group.\qquad $\square$
\section{The EDS associated with the (2,3,5) distribution associated with the movement of three ants}   \label{se5}
In this section we again analyze the movement of the three ants according to the {\bf rule B}, but now from the Exterior Differential System (EDS) point of view. Our goal here will be to give proofs of Theorem \ref{33} and Proposition \ref{pro}, which will be totally independent of the geometric arguments in Section \ref{sec4}. In particular, we aim to charactarize the $(2,3,5)$ distribution $\mathcal{D}$ associated with ants moving according to {\bf Rule B} in terms of the corresponding EDS, and use it to \emph{calculate} the formula for the Cartan quartic of $\mathcal{D}$. One can try to make this calculations in terms of the original Euclidean coordinates $(x_1,y_1,x_2,y_2,x_3,y_3)$ on the plane, but we doubt that it is possible in any finite time (even with the help of very powerful computer equipped with a good symbolic calculation program). It is why here we use totally new parametrization of the ants system, using the obsevation that its symmetry group is at least as large as the semidirect product $\slg(2,\bbR)\rtimes \bbR^2$ of the groups $\slg(2,\bbR)$ and $\bbR^2$. Actually we embed the non-controllable 6-dimensional configuration space of the ants moving according to {\bf rule B} in the 6-dimensional group $\glg(2,\bbR)\rtimes \bbR^2$ and identify a chosen 5-dimensional orbit $M_A$ of controllable movement with the group $\slg(2,\bbR)\rtimes \bbR^2$. The advantage of such an approach is that given the identification of the configuration space manifold $M_A$ with a Lie group manifold, we have a preferred set of 1-forms on $M_A$. This is given by the \emph{coframe} of \emph{left invariant forms}, known as the \emph{Maurer-Cartan forms}. These are \emph{God given} on any Lie group (see \cite{helgason} Ch.II $\S$ 7, 135-138). In the coordinates on $M_A$ adapted to the group parameters they are easy to calculate explicitely. Having them in coordinates adapted to the symmetry of the system will be extremely helpful to find the $(2,3,5)$ geometric inavriant 1-forms constituting the EDS  described in \cite{nan2}, p. 94, encapsulating the curvature and Cartan quartic of the ants-under-{\bf Rule}-{\bf B} $(2,3,5)$ distribution $\mathcal D$. Having this EDS we will be able to determine the Cartan quartic of $\mathcal D$ and its type.

\subsection{Parametrization in terms of $\glg(2,\bbR)\rtimes \bbR^2$} Let us make use of the fact that the symmetry of the ant system obeying {\bf rule B} is at least $\slg(2,\bbR)\rtimes \bbR^2$ and that the reduced system can be considered on a manifold $M_A$ being diffeomorphic to this group. Actually we paramatrize all triangles on the plane by the elements of the 6-dimensional $\glg(2,\bbR)\rtimes \bbR^2$ group.

For this we take a `standard triangle'\footnote{Note that \emph{any} choice of a `standard triangle' will work; we have chosen this one because we like the origin of the Cartesian coordinate system, and the isosceles rectangular triangles.} in $\bbR^2$, which we define in terms of its vertices at $\vec{r}_1=(0,0)$, $\vec{r}_2=(1,0)$ and $\vec{r}_3=(0,1)$ in a chosen Cartesian system $(x,y)$ on the plane. Now, any other triangle on the plane is obtained by acting on the vertices of the standard triangle with the group $\glg(2,\bbR)\rtimes \bbR^2$. It is convenient to represent the group $\glg(2,\bbR)\rtimes \bbR^2$ as a group of invertible $3\times 3$ real matrices
\be h=\bma a&b&x\\p&q&y\\0&0&1\ema,\quad aq-bp\neq 0,\label{det}\ee
and vectors $\vec{r}$ in $\bbR^2$ as column vectors $r=(\vec{r},1)^t$. Then the action of the group $\glg(2,\bbR)\rtimes \bbR^2$ in $\bbR^2$ can be read off from the action $(h,r)\to h.r$ coming from the multiplication `.' of matrices $h$ from $\glg(2,\bbR)\rtimes \bbR^2$ and vectors $r$ of the form $r=(\vec{r},1)^t$ from $\bbR^3$. Using this action, we can transform the standard triangle with vertices $\vec{r}_i$, $i=1,2,3$, to any other triangle on the plane. In this way the most general triangle on the plane will have the vertices $$\begin{aligned}
  \vec{R}_1=h.r_1&=(x,y),\quad\quad \vec{R}_2=h.r_2=(a+x,p+y),\\&\vec{R}_3=h.r_3=(b+x,q+y).\end{aligned}$$

Thus we locally parametrized the triangles on the plane by coordinates $(x,y,a,b,p,q)$ in an open set of $\bbR^6$ identified with the affine transformation group of the plane  $\glg(2,\bbR)\rtimes \bbR^2$. We put the three ants each one in one of the vertices $\vec{R}_1,\vec{R}_2,\vec{R}_3$ of the general triangle. We are interested in curves $m(t)=(x(t),y(t),a(t),b(t),p(t),q(t))$ in $\glg(2,\bbR)\rtimes \bbR^2$, corresponding to the movement of the ants in time, and we want that these curves satisfy:
$$\frac{\der \vec{R}_1}{\der t}\, ||\, (\vec{R}_2-\vec{R}_3),\quad\frac{\der \vec{R}_3}{\der t}\, ||\,(\vec{R}_1-\vec{R}_2),\quad\frac{\der \vec{R}_2}{\der t}\, ||\,(\vec{R}_3-\vec{R}_1),$$
which is the implementation of {\bf rule B}. More explicitly this rule means that
\be \begin{aligned}(\dot{x},\dot{y})\, ||\,(a-b,&p-q),\quad(\dot{a}+\dot{x},\dot{p}+\dot{y})\, ||\,(b,q),\\&(\dot{b}+\dot{x},\dot{q}+\dot{y})\, ||\,(-a,-p).\end{aligned} \label{mv}\ee

Instead of writing the $(3,6)$ velocity distribution $\Db$ of the ants as in Section \ref{B}, we now write this distribution in terms on three vector fields on  $\glg(2,\bbR)\rtimes \bbR^2$ corresponding to the three \emph{primitive moves}.

Obviously the following three moves satisfy {\bf rule B}:
\begin{itemize}
\item[]{\bf Move 1:} The first ant, with a position at $\vec{R}_1$ moves, the other two ants are at rest,
\item[]{\bf Move 2:} The second ant, with a position at $\vec{R}_2$ moves, the other two ants are at rest,
  \item[]{\bf Move 3:} The third ant, with a position at $\vec{R}_3$ moves, the other two ants are at rest.
\end{itemize}
Looking at the implementation of {\bf rule B} in \eqref{mv} we see that:\\
\noindent
$\bullet$ {\bf Move 1} means that $\dot{x}=c(a-b)$, $\dot{y}=c(p-q)$, $\dot{a}=-\dot{x}$, $\dot{p}=-\dot{y}$,$\dot{b}=-\dot{x}$,$\dot{q}=-\dot{y}$, resulting in $(\dot{x},\dot{y},\dot{a},\dot{b},\dot{p},\dot{q})=c(a-b,p-q,b-a,b-a,q-p,q-p)$, or a generating vector field \be V_1=(a-b)\big(\partial_x-\partial_a-\partial_b\big)+(p-q)\big(\partial_y-\partial_p-\partial_q\big).\label{v1}\ee
\noindent
 $\bullet$ {\bf Move 2} means that $\dot{x}=\dot{y}=\dot{b}=\dot{q}=0$, $\dot{a}=c b $, $\dot{p}=c q$, resulting in $(\dot{x},\dot{y},\dot{a},\dot{b},\dot{p},\dot{q})=c(0,0,b,0,q,0)$, or a generating vector field \be V_2=b \partial_a+q\partial_p.\label{v2}\ee
    \noindent
$\bullet$ {\bf Move 3} means that $\dot{x}=\dot{y}=\dot{a}=\dot{p}=0$, $\dot{b}=-c a $, $\dot{q}=-c p$, resulting in $(\dot{x},\dot{y},\dot{a},\dot{b},\dot{p},\dot{q})=c_1(0,0,0,-a,0,-p)$, or a generating vector field \be V_3=-a \partial_b-p\partial_q.\label{v3}\ee
The general move of the ants according to {\bf rule B} is a superposition of these three primitive moves, so the velocity distribution $\Db$ of the ants moving in $\glg(2,\bbR)\rtimes \bbR^2$ according to {\bf rule B} is
$$\Db=\Span(V_1,V_2,V_3),$$
where the vector fields $V_1,V_2,V_3$ are given by the respective formulas \eqref{v1}, \eqref{v2} and \eqref{v3}.

Now, let us consider the determinant $$Det(h)=a q- bp$$
  of the $\glg(2,\bbR)\rtimes \bbR^2$ valued matrix $h$ as in \eqref{det}. One can easilly check that a general vector field
  $$V=f_1 V_1+f_2 V_2+f_3 V_3$$
  from the distribution $\Db$ annihilates $Det(h)$, i.e.\footnote{Note that in the formula ${\mathcal L}_V\big(Det(h)\big)=0$, we used the explicit forms of vector fields $V_1$, $V_2$ and $V_3$ which are given explicitely as \emph{differential operators} in formulas \eqref{v1}, \eqref{v2} and \eqref{v3}. Thus the expression $V\big(Det(h)\big)$ means action of the first order differential operator $V$ on a \emph{function} $Det(h)=aq-bp$. Direct simple calculation using the explicit form of $V_1$, $V_2$ and $V_3$ shows that ${\mathcal L}_V\big(Det(h)\big)=0$.}
  $${\mathcal L}_V\big(Det(h)\big)=0,\quad \forall V\in\Db.$$
  Here ${\mathcal L}_V\big(f\big)$ denotes the Lie derivative of a function $f$ with respect to the vector field $V$.

  Thus, as we already know, the distribution $\Db$ is tangent to 5-dimensional submanifolds in $\glg(2,\bbR)\rtimes \bbR^2$ consisting of the elements with \emph{constant determinant}. Each of these submanifolds is diffeomorphic to the $\slg(2,\bbR)\rtimes \bbR^2$ group, i.e. the group of motions in $\bbR^2$ preserving volumes.

  Let us now introduce a foliation of $\glg(2,\bbR)\rtimes \bbR^2$ by 5-dimensional submanifolds $N_s$ defined by
  $$
  \resizebox{0.49\textwidth}{!}{$N_s=\{h\,\,{\mathrm as}\,\,\mathrm{ in}\,\, \eqref{det}\,\,\mathrm{s.t.}\,\,Det(h)=s=\mathrm{const}\}\subset \glg(2,\bbR)\rtimes \bbR^2.$}$$
On each submanifod $N_s$ we then have a rank 2-distribution
  $$\Da=\Span(V_1-V_3,V_2-V_1),$$
such that it is the square root of $\Db$, $\Db=[\Da,\Da]+\Da$. The distribution $\Da$ obviously has the growth vector $(2,3,5)$ on each 5-dimensional manifold $N_s$.

  \subsection{The $(2,3,5)$ distribution and the Maurer-Cartan forms on $\glg(2,\bbR)\rtimes \bbR^2$} Since the distribution $\Da$ is defined on the 5-dim submanifolds of $\glg(2,\bbR)\rtimes \bbR^2$ diffeomorphic to $\slg(2,\bbR)\rtimes \bbR^2$, and since this distribution is definitely $\slg(2,\bbR)\rtimes \bbR^2$ invariant, it is natural to ask how to define it in terms of the Maurer-Cartan forms on $\glg(2,\bbR)\rtimes \bbR^2$ or $\slg(2,\bbR)\rtimes \bbR^2$.

  If we view of $\glg(2,\bbR)\rtimes \bbR^2$ as the group of matrices $h$ defined in \eqref{det}, then the basis $(\tau^1,\tau^2,\tau^3,\tau^4,\tau^5,\tau^6)$ of the Maurer-Cartan forms on $\glg(2,\bbR)\rtimes \bbR^2$ can be easilly found from the formula
  $$h^{-1}\der h= \sum_{i=1}^6\tau^i E_i,$$
  with
  \be \begin{aligned}
    E_1=&\bma 0&1&0\\0&0&0\\0&0&0\ema,\quad E_2=\bma 0&0&0\\1&0&0\\0&0&0\ema,\quad E_3=\bma 1&0&0\\0&-1&0\\0&0&0\ema,\\E_4=&\bma 0&0&1\\0&0&0\\0&0&0\ema,\quad E_5=\bma 0&0&0\\0&0&1\\0&0&0\ema,\quad E_6=\bma 1&0&0\\0&1&0\\0&0&0\ema.\end{aligned}\label{maci}\ee
  Explicitly we have
  $$\begin{aligned}
    &\tau^1=\frac{q\der b-b\der q}{aq-bp},\quad\tau^2=\frac{a\der p-p\der a}{aq-bp},\\&\tau^3=\frac{q\der a-a\der q+p\der b-b\der p}{2(aq-bp)},\quad\tau^4=\frac{q\der x-b\der y}{aq-bp},\\&\tau^5=\frac{a\der y-p\der x}{aq-bp},\quad \tau^6=\frac{\der(aq-bp)}{2(aq-bp)}.
  \end{aligned}
  $$
Note that this basis of the Maurer-Cartan forms has the property that on each of the 5-dimensional submanifolds $N_s$, where $s=aq-bp=\mathrm{const}$ the sixth 1-form $\tau^6$ identically vanish, $$\tau^6\equiv 0\quad\mathrm{on}\,\,\mathrm{each}\quad N_s.$$
Now we can look for the most general linear combination (with constant coefficients!) of the Maurer-Cartan forms $(\tau^1,\tau^2,\tau^3,\tau^4,\tau^5)$  which \emph{annihilates} the $(2,3,5)$ distribution $\Da$ on each $N_s$. Since the explicit form of the vector fields $V_1-V_3$ and $V_2-V_1$ spanning the $(2,3,5)$ distribution $\Da$ in coordinates $(a,b,p,q,x,y)$ reads:
$$\resizebox{0.47\textwidth}{!}{$\begin{aligned}&V_1-V_3=(a-b)(\partial_x-\partial_a)+(p-q)(\partial_y-\partial_p)+b\partial_b+q\partial_q,\\& V_2-V_1=(a-b)(\partial_b-\partial_x)+(p-q)(\partial_q-\partial_y)+a\partial_a+p\partial_p,\end{aligned}$}$$
and since these vector fields annihilate $\tau^6$, one easily finds that the basis of the annihilator of $\Da$ on $N_s$ is given by
$$\Da^\perp=\Span(\tau^3-\tau^5,\tau^4+\tau^5,\tau^1-\tau^2-\tau^3).$$
Of course, since we restrict our attention to the leaves $N_s$ of the foliation of $M$, the Maurer-Cartan forms appearing in this formula have $aq-bp=s=const$.
This means that we can define the $(2,3,5)$ distribution $\Da$ on each $N_s$ as the annihilator of the forms
$$\resizebox{0.5\textwidth}{!}{$\begin{aligned}
   s& \theta^1=s(\tau^3-\tau^5)=\,p\der x-a\der y+\tfrac12(q\der a+p\der b-b\der p-a\der q),\\
 s &\theta^2=s(\tau^4+\tau^5)=\,-(p-q)\der x+(a-b)\der y,\\
2&s\theta^3=2s(\tau^1-\tau^2-\tau^3)=\\&\quad\quad(2p-q)\der a+(2q-p)\der b-(2a-b)\der p-(2b-a)\der q,\\
\end{aligned}$}$$
These three 1-forms can be supplemented by the Maurer-Cartan forms
$$\begin{aligned}
  s\theta^4=-\tau^2=&\,p\der a-a\der p,\\
  s\theta^5=\tau^2-\tau^1=&\,-p\der a-q\der b+a\der p+b\der q,
\end{aligned}
$$
to a coframe $(\theta^1,\theta^2,\theta^3,\theta^4,\theta^5)$ on $N_s$. In this coframe the distribution
$$\Da=\{V\in \Gamma({\mathrm T}N_s)\,:\, V\hook\theta^1=V\hook\theta^2=V\hook \theta^3=0)\}$$
and because of the boxed terms in the formulas below
\be\resizebox{0.47\textwidth}{!}{$\begin{aligned}
\der\theta^1&=\theta^1\dz(\theta^3+\theta^4+\theta^5)+\theta^2\dz\theta^4+\boxed{\theta^3\dz\theta^4},\\
\der\theta^2&=-\theta^1\dz(2\theta^3+\theta^5)-\theta^2\dz(\theta^3+\theta^4+\theta^5)+\boxed{\theta^3\dz\theta^5},\\
\der\theta^3&=-\theta^3\dz(4\theta^4+2\theta^5)+\boxed{3\theta^4\dz\theta^5},\\
\der\theta^4&=-2\theta^3\dz\theta^4+2\theta^4\dz\theta^5,\\
\der\theta^5&=\theta^3\dz(4\theta^4+2\theta^5)-4\theta^4\dz\theta^5,\label{sys}
\end{aligned}$}
\ee
its growth vector is visibly\footnote{The appearence of the boxed terms with \emph{nonzero constant} coefficients implies - by the Maurer-Cartan formula (see \cite{helgason}  Ch.II $\S$ 7, Proposition 7.2, p. 137) that the vector fields $X_i$ dual to the forms $\theta^i$, $X_i\hook\theta^j=\delta_i{}^j$, satisfy in particular: $[X_4,X_5]=-3 X_3\mod X_4,X_5$, $[X_3,X_5]=- X_2\mod X_3,X_4,X_5$ and $[X_3,X_4]=- X_1\mod X_3,X_4,X_5$. This shows that the rank 2 distribution spanned by the vector fields $X_4$ and $X_5$ is precisely $(2,3,5)$ and can be identified with the ants distribution $\mathcal D$.} $(2,3,5)$. Also, the appearence of only constant coefficients in \eqref{sys}, visibly shows that the distribution $\Da$ is \emph{homogeneous} with the symmetry group being at least as large as $\slg(2,\bbR)\rtimes\bbR^2$.

This proves the following proposition.

\begin{proposition}
  The $(2,3,5)$ distribution $\Da$ of three ants moving on the floor according to {\bf rule B} is locally equivalent to a homogeneous $(2,3,5)$ distribution on the Lie group  $\slg(2,\bbR)\rtimes\bbR^2$, which is defined as the annihilator of the three Maurer-Cartan forms $(\theta^1,\theta^2,\theta^3)$ with  $$\theta^1=\tau^3-\tau^5, \quad\theta^2=\tau^4+\tau^5, \quad\theta^3=\tau^1-\tau^2-\tau^3,$$ where $\tau^i$s are defined in terms of the general element $h\in \slg(2,\bbR)\rtimes\bbR^2$, as in \eqref{det}, and the basis $E_i$ in $\gla(2,\bbR)\roplus \bbR^2$, as in \eqref{maci}, by $h^{-1}\der h=\sum_{i=1}^6\tau^iE_i$.
\end{proposition}

\subsection{Cartan quartic for the ants' $(2,3,5)$ distribution on $\slg(2,\bbR)\rtimes\bbR^2$}
The Cartan quartic \cite{cartan} for the ants' distribution $\Da$ can be computed in various ways. Here we do it by calculating explicitly the conformal $(2,3)$ signature metric \cite{nurdif} associated with $\Da$. It follows, that in the coframe $(\theta^1,\theta^2,\theta^3,\theta^4,\theta^5)$ satisfying the system \eqref{sys}, the conformal representative of this metric can be taken as:
\be\begin{aligned}
g=&\theta^1\,\big(\,\,45\,\theta^5+60\,\theta^3+27\,\theta^2+27\,\theta^1\,\,\big)\,-\\&\,\theta^2\,\big(\,\,45\,\theta^4-30\,\theta^3-27\,\theta^2\,\,\big)\,+\,10\,(\,\theta^3\,)^2,\end{aligned}
\label{metric}
\ee
with the product between the 1-forms above being $\lambda\mu=\tfrac12(\lambda\otimes\mu+\mu\otimes\lambda)$.

Calculating the Weyl tensor\footnote{Here we use the basic invariants of the conformal pseudo-riemannian geometry, namely the \emph{Weyl tensor} and the \emph{conformally Einstein condition}. The conformal geometry and the theory of its invariants is a big subject on its own. A convenient shortcut to the definitions of conformal objects and notions we use here, can be found in Sections 2.1-2.2 of \cite{gover}.} of this metric in the null coframe, and using the procedure of calculating the Cartan quartic from the conformal metric described in \cite{nan}, we find that the Cartan quartic of the corresponding to $(2,3,5)$ distribution $\Da$ is \emph{of type} $D$ in the parabolic geometric language \cite{franco} or, what is the same, \emph{has no real roots}. Morevoer the metric is \emph{not} conformal to an Einstein metric. This in particular means that the distribution $\Da$ is NOT $G_2$ flat\footnote{There are many meanings of the word \emph{flatness} in mathematics. Usually it refers to the \emph{zero curvature of a} certain \emph{connection}. It is the meaning we use in this sentence: The Cartan $\mathfrak{g}_2^*$-valued Cartan connection of the $(2,3,5)$ distribution $\mathcal{D}$  associated with ants moving on a 5-dimensional orbit according to {\bf rule B} has non-zero curvature, i.e. it is \emph{not flat}. This has nothing to do with the word `flatness' used in the nonlinear control theory.}. Also, using the \emph{Cartan reduction procedure} for the Cartan system associated to the distribution $\Da$, as in Theorem 8 in \cite{nurdif}, we established that the districbution $\Da$ has precisely 5-dimensional Lie algebra of symmetries. This proves Theorem \ref{33} from the end of Section \ref{B}.
\section{Acknowledgements}
The idea of this work arose during discussions between the second author and Gil Bor at his house in Guanajuato, MX, and also at his office in CIMAT, in November 2015. PN is very grateful to Gil Bor and all the members of the CIMAT institute in Gaunajuato for creating a very friendly working atmosphere and partial financial cover of his stay while in Mexico.

\end{document}